\documentclass[10pt,a4paper]{amsart}

\usepackage{amssymb}
\usepackage{graphicx}
\usepackage[ps2pdf]{hyperref} 
\hypersetup{colorlinks=true, linkcolor=blue, anchorcolor=blue, 
citecolor=red, filecolor=blue, menucolor=blue,urlcolor=blue}
\numberwithin{equation}{section}
\newtheorem{theorem}{Theorem}[section]
\newtheorem{corollary}[theorem]{Corollary}
\newtheorem{proposition}[theorem]{Proposition}

\theoremstyle{remark}
\newtheorem{remark}{Remark}[section]

\theoremstyle{definition}

\newtheorem{example}[theorem]{Example}

\begin{document}

\title%
[Product manifolds]%
{Dispersive estimates and NLS on product manifolds}

\begin{abstract}
We prove a general dispersive estimate for a Schr\"odinger
type equation on a product manifold,
under the assumption that the equation restricted to each factor
satisfies suitable dispersive estimates. Among the
applications are the two-particle
Schr\"odinger equations
\begin{equation*}
  iu_{t}-\Delta_{x,y}u+V(x-y)u=0
\end{equation*}
on $\mathbb{R}^{2n}$,
and the nonlinear Schr\"odinger equation on the product
of two real hyperbolic spaces $\mathbb{H}^{m}\times\mathbb{H}^{n}$.
\end{abstract}

\date{\today}

\author{Vittoria Pierfelice}
\address{ Laboratoire MAPMO UMR 6628, 
Universit\'e d'Orl\'eans, 
Batiment de math\'ematiques - Route de Chartres
B.P. 6759 - 45067 Orl\'eans cedex 2
FRANCE}
\email{vittoria.pierfelice@univ-orleans.fr}



\maketitle

\section{Introduction}\label{sec:intro}

Let $X=M \times N$ be a product of oriented riemannian manifolds, each
endowed with its canonical volume form, and let us consider
three unbounded selfajoint operators, $L$ on $L^{2}(X)$, 
$H$ on $L^{2}(M)$ and $K$ on $L^{2}(N)$. We shall assume that the 
operator $L$ is the \emph{sum} of $H$ and $K$, in the following
sense: we assume that, for every couple
of functions $f,g$, with $f(x)$ in a dense subset of the domain of $H$ and 
$g(y)$ in a dense subset of the domain of $K$, we have that
the function $f(x)g(y)$ is in the domain of $L$ and 
\begin{equation}\label{eq:basicprop}
  L (f(x)g(y))=Hf \cdot g+f \cdot Kg.
\end{equation}
This situation is quite common and occurs in a number of interesting
and natural examples. We mention a few:

\begin{example}\label{exa:1}
  The simplest case of course is given by the standard Laplacian
  on $\mathbb{R}^{n}_{x}\times \mathbb{R}^{m}_{y}$;
  we have
  \begin{equation*}
    \Delta_{x,y}=\Delta_{x}+\Delta_{y}
  \end{equation*}
  and the relation \eqref{eq:basicprop} is satisfied with
  the choices $L=-\Delta_{x,y}$, $H=-\Delta_{x}$, $K=-\Delta_{y}$.
  In greater generality, we can choose
  $L$, $H$ and $K$ to be the Laplace-Beltrami operators on 
  the three manifolds $X,M$ and $N$
  respectively. Indeed, in local coordinates the metric on $X$
  is given by a block matrix, with two blocks corresponding to the
  metrics of $M$ and $N$; using the explicit representation of the
  Laplace-Beltrami operators it is easy to check that
  \begin{equation*}
    \Delta_{X}(f(x)g(y))=
    \Delta_{M}f \cdot g+f \cdot\Delta_{N}g.
  \end{equation*}
\end{example}

\begin{example}\label{exa:2}
  On $\mathbb{R}^{m}_{x}\times \mathbb{R}^{n}_{y}$,
  consider the Schr\"odinger operator
  \begin{equation}\label{eq:splitpot}
    L= -\Delta_{x,y}+U(x,y),\qquad U(x,y)=
    V(x)+W(y)
  \end{equation}
  where the potential $U(x,y)$ can be split in the sum of
  two potentials depending only on
  a group of variables each. Then we may choose
  \begin{equation*}
    H=-\Delta_{x}+V(x),\qquad K=-\Delta_{y}+W(y).
  \end{equation*}
  More generally, $H$ and $K$ can be two \emph{electromagnetic}
  Schr\"odinger operators of the form
  \begin{equation*}
    (i \nabla_{x}-A(x))^{2}+V(x),\qquad
    (i \nabla_{y}-B(y))^{2}+W(y)
  \end{equation*}
  with $A:\mathbb{R}^{m}\to \mathbb{R}^{m}$ and
  $B:\mathbb{R}^{n}\to \mathbb{R}^{n}$.
\end{example}

\begin{example}\label{exa:3}
  The wave function $u(t,x,y)$
  of two interacting particles is governed by 
  a Schr\"odinger equation of the form
  \begin{equation}\label{eq:2part}
    iu_{t}-\Delta_{x,y}u+V(x-y)u=0,\qquad t\in \mathbb{R},\quad
      (x,y)\in \mathbb{R}^{3+3}.
  \end{equation}
  By the change of variables $x'=x+y$, $y'=x-y$, equation
  \eqref{eq:2part} reduces to the following equation for
  $v(t,x',y')=u(t,x,y):$
  \begin{equation*}
    iv_{t}-\Delta_{x',y'}v+V(y')v=0.
  \end{equation*}
  We see that the Schr\"odinger operator here belongs to the
  class considered in Example \ref{exa:2}.
\end{example}

The first goal of this paper is to show by an elementary
abstract argument that the dispersive properties
of the flows $e^{itL}$, $e^{itH}$ and $e^{itK}$ are related in a
natural way. This approach allows to handle some cases when the
usual methods to prove dispersive estimates can not be
applied. Although the methods are completely elementary,
the result has a number of interesting consequences.
Our basic result is the following:

\begin{theorem}\label{the:1}
  Assume the Schr\"odinger flows for $H$ and $K$ satisfy,
  for some real $a,b\ge0$, and for $t$ belonging
  to an interval $I \subseteq \mathbb{R}$,
  dispersive estimates of the form
  \begin{equation}\label{eq:assdis}
    \|e^{itH}\phi\|_{L^{r}(M)}\lesssim
       |t|^{-a}\|\phi\|_{L^{\widetilde{r}}(M)},\qquad
    \|e^{itK}\psi\|_{L^{r}(N)}\lesssim
       |t|^{-b}\|\psi\|_{L^{\widetilde{r}}(N)}.
  \end{equation}
  for some exponents $\widetilde{r}\le r$ in $[1,\infty]$.
  Then the flow of $L$ satisfies for $t\in I$ the estimate
  \begin{equation}\label{eq:Ldis}
    \|e^{itL}f\|_{L^{r}(M \times N)}\lesssim
       |t|^{-a-b}\|f\|_{L^{\widetilde{r}}(M \times N)}.
  \end{equation}
\end{theorem}

It is always possible to interpolate the previous
dispersive estimate with the conservation of energy
\begin{equation*}
  \|e^{itL}f\|_{L^{2}}\equiv\|f\|_{L^{2}},
\end{equation*}
which follows from the selfadjointness of $L$.
In particular, if the assumptions of Theorem \ref{the:1}
hold with $r=\infty$, $\widetilde{r}=1$,
we obtain the
complete set of dispersive $L^{q'}-L^{q}$ estimates
\begin{equation}\label{eq:LqLqp}
  \|e^{itL}f\|_{L^{q}(M \times N)}\lesssim
     |t|^{-(a+b)\left(1- \frac2q\right)}
     \|f\|_{L^{q'}(M \times N)},\qquad
     2\le q\le \infty.
\end{equation}
Following the methods of
\cite{Kato87-a},
\cite{GinibreVelo95-a},
\cite{KeelTao98-a}, it is then possible
to deduce in a standard way the 
corresponding \emph{Strichartz estimates}. We use the notation,
for any finite or infinite interval $I \subseteq \mathbb{R}$,
\begin{equation*}
  \|F(t,x,y)\|_{L^{p}_{I}L^{q}_{\phantom I}}=
  \left(
  \int_{I}
  \left(
    \int_{M \times N}|F(t,x,y)|^{q}dV_{x,y}
  \right)^{\frac pq}dt
  \right)^{\frac1p}.
\end{equation*}
We also define an \emph{admissible couple}, associated to the
index $a+b$, as follows: when $a+b>1$,
the couple $(p,q)$ is admissible if
it satisfies the conditions
\begin{equation}\label{eq:adm}
  \frac1p+\frac{a+b}{q}=\frac{a+b}{2},\qquad
  2\le p\le \infty,\qquad
  \frac{2(a+b)}{a+b-1}\ge q\ge 2;
\end{equation}
when $0<a+b\le1$, the conditions are
\begin{equation}\label{eq:adm2}
  \frac1p+\frac{a+b}{q}=\frac{a+b}{2},\qquad
  \frac{2}{a+b}< p\le \infty,\qquad
  \infty> q\ge 2.
\end{equation}
We also denote with $q'$ the dual exponent to $q$. 
The value
\begin{equation*}
  (p,q)=\left(2,
    \frac{2(a+b)}{a+b-1}
  \right)
\end{equation*}
(when $a+b\ge1$)
is the \emph{endpoint}; notice that $q\neq \infty$
in all cases considered here.

Then we have:

\begin{proposition}\label{pro:stri}
  Assume $X,M,N$ and $L,H,K$ are as in Theorem \ref{the:1}
  with $r=\infty$, $\widetilde{r}=1$. Then
  the following estimates hold:
  \begin{equation}\label{eq:homstr}
    \left\|
      e^{itL}f
    \right\|_{L^{p}_{I}L^{q}_{\phantom I}}
    \lesssim
    \|f\|_{L^{2}(M \times N)},
  \end{equation}
  \begin{equation}\label{eq:nonhomstr}
    \left\|
      \int_{0}^{t}e^{i(t-s)L}F(s,x,y)ds
    \right\|_{L^{p}_{I}L^{q}_{\phantom I}}
    \lesssim
    \|F\|_{L^{\widetilde{p}'}_{I}L^{\widetilde{q}'}_{\phantom I}}
  \end{equation}
  for all admissible couples $(p,q)$ and $(\widetilde{p},\widetilde{q})$.
\end{proposition}


\begin{remark}\label{rem:kterms}
  It is evident that the result generalizes to a finite product
  of manifolds $X=M_{1}\times \dots \times M_{k}$ and
  an operator $L$ on $X$
  decomposable as the sum $L=H_{1}+\dots+H_{k}$
  where each $H_{j}$ acts on $M_{j}$ only.
\end{remark}

\begin{remark}\label{rem:deriv}
  If the flows satisfy estimates with loss of derivatives
  of the form
  \begin{equation*}
    \|e^{itH}\phi\|_{L^{\infty}(M)}\lesssim
       |t|^{-a}\|H^{r}\phi\|_{L^{1}(M)},\qquad
    \|e^{itK}\psi\|_{L^{\infty}(N)}\lesssim
       |t|^{-b}\|H^{r}\psi\|_{L^{1}(N)}
  \end{equation*}
  then it is easy to extend the result of Theorem \ref{the:1}
  and obtain the estimate
  \begin{equation*}
    \|e^{itL}f\|_{L^{\infty}(M \times N)}\lesssim
       |t|^{-a-b}\|H^{r}K^{s}f\|_{L^{1}(M \times N)}.
  \end{equation*}
  To this end, it is sufficient to apply the argument in the proof
  to the modified flows
  \begin{equation*}
    H^{-r}e^{itH} \quad\text{and}\quad
    K^{-s}e^{itK}
  \end{equation*}
  instead of $e^{itH},e^{itK}$.
\end{remark}

Despite its simplicity,
Theorem \ref{the:1} has several applications. 
We begin by studying the case of Schr\"odinger operators with
potential perturbations on $\mathbb{R}^{n}$. A first example
is based on the 1D decay results of 
\cite{DanconaFanelli06-a}, \cite{GoldbergSchlag04-b}:

\begin{corollary}\label{cor:3}
  Let $V(x)\ge0$ be a real valued function such that
  \begin{equation}\label{eq:weight}
    (1+|x|)^{2}V(x)\in L^{1}(\mathbb{R}).
  \end{equation}
  Then, for all $n\ge1$,
  the solution of the Schr\"odinger equation on $\mathbb{R}^{n}$
  \begin{equation*}
    iu_{t}-\Delta u+(V(x_{1})+\dots+V(x_{n}))u=0,\qquad u(0,x)=f(x)
  \end{equation*}
  satisfies the estimate
  \begin{equation*}
    |u(t,x)|\lesssim |t|^{-n/2}\|f\|_{L^{1}(\mathbb{R}^{n})}.
  \end{equation*}
\end{corollary}

We notice that in dimension $n=2$ this gives a classes of
potentials for which a sharp $L^{1}-L^{\infty}$ estimate
is true; no other classes are known to our knowledge (the
only known estimates are of type $L^{p}-L^{p'}$ with
$2\le p<\infty$, see \cite{Yajima99-a}).

In dimension $n\ge3$, it is not known what are the
optimal conditions on a potential $V(x)$ such that 
the flow $i \partial_{t}-\Delta+V(x)$ satisfies a dispersive
estimate. However there are several sufficient conditions due to
different authors. We mention for instance the following
(see \cite{Yajima95-a}):
$n\ge3$, $p_{0}>n/2$, $\delta>3n/2+1$,
$\ell_{0}=0$ if $n=3$ and $\ell_{0}=[(n-1)/2]$ if $n\ge4$, and
$V:\mathbb{R}^{n}\to \mathbb{R}$ satisfies
\begin{equation}\label{eq:yaj}
  \|D^{\alpha}V\|_{L^{p_{0}}(|x-y|\le1)}\le 
  \frac{C}{(1+|x|)^{\delta}}\qquad
  \forall|\alpha|\le \ell_{0};
\end{equation}

\begin{corollary}\label{cor:4}
  Let $m,n\ge3$, assume the potentials $V:\mathbb{R}^{m}\to \mathbb{R}$
  and $W:\mathbb{R}^{n}\to \mathbb{R}$ satisfy condition \eqref{eq:yaj}
  (in dimension $m$ and $n$ respectively). Then the solution $u(t,x,y)$
  of the Schr\"odinger equation on $\mathbb{R}^{m+n}$
  \begin{equation*}
    iu_{t}-\Delta_{x,y}u+V(x)u+W(y)u=0,\qquad
    u(0,x,y)=f(x,y)
  \end{equation*}
  satisfies the dispersive estimate
  \begin{equation*}
    |u(t,x,y)|\lesssim|t|^{-\frac{m+n}{2}}\|f\|_{L^{1}(\mathbb{R}^{m+n})}.
  \end{equation*}
\end{corollary}

For the Schr\"odinger equation describing the interaction
of two particles, which was examined in Example \ref{exa:3},
we can prove the following:

\begin{corollary}\label{cor:5}
  Let $n\ge3$ and $V{:\mathbb{R}^{n}\to \mathbb{R}}$ satisfying
  condition \eqref{eq:yaj}. Then the solution
  of the Schr\"odinger equation on 
  $\mathbb{R}^{2n}=\mathbb{R}^{n}_{x}\times \mathbb{R}^{n}_{y}$
  \begin{equation*}
    iu_{t}-\Delta_{x,y}u+V(x-y)u=0,\qquad
    u(0,x,y)=f(x,y)
  \end{equation*}
  satisfies the dispersive estimate
  \begin{equation*}
    |u(t,x,y)|\lesssim|t|^{-n}\|f\|_{L^{1}(\mathbb{R}^{2n})}.
  \end{equation*}
\end{corollary}

As a final application, we consider a manifold $X$ which is the product
of two hyperbolic spaces
\begin{equation*}
  X=\mathbb{H}^{m}\times \mathbb{H}^{n},\qquad
  m,n\ge2.
\end{equation*}
The Schr\"odinger equation on hyperbolic spaces was investigated
in several papers; in particular, weighted Strichartz estimates were
proved in \cite{Pierfelice06-b} while sharp dispersive estimates were
obtained in \cite{AnkerPierfelice09-a}.  We recall the main
result of \cite{AnkerPierfelice09-a} for 
$e^{it \Delta_{\mathbb{H}^{n}}}$: for all 
$r,\widetilde{r}\in(2,\infty]$ we have
\begin{equation}\label{eq:AP}
  \|e^{it \Delta_{\mathbb{H}^{n}}}f\|_{L^{r}(\mathbb{H}^{n})}\lesssim
  \begin{cases}
    |t|^{-\max\{\frac12-\frac1r,\frac12-\frac1{\widetilde{r}}\}n} 
    \|f\|_{L^{\widetilde{r}'}}
    &\text{if $ 0<|t|\le1 $,}\\
    |t|^{-\frac32}\|f\|_{L^{\widetilde{r}'}} &\text{if $ |t|\ge1 $.}
  \end{cases}
\end{equation}
Using Theorem \ref{the:1} we obtain:

\begin{corollary}\label{cor:6}
  Consider the Schr\"odinger equation
  $iu_{t}-\Delta_{X} u=0$, $u(0)=f$ on the product manifold
  $X=\mathbb{H}^{m}\times \mathbb{H}^{n}$, $m,n\ge2$, where
  $\Delta_{X}$ is the Laplace-Beltrami operator on $X$.
  Then the solution $u(t,x,y)$ satisfies, for all 
  $r,\widetilde{r}\in(2,\infty]$, the dispersive estimate
  \begin{equation}\label{eq:prodman}
    \|u(t)\|_{L^r}\lesssim
    \begin{cases}
      |t|^{-\max\{\frac12-\frac1r,\frac12-\frac1{\widetilde{r}}\}(n+m)} 
      \|f\|_{L^{\widetilde{r}'}}
      &\text{if $ 0<|t|\le1 $,}\\
      |t|^{-3}\|f\|_{L^{\widetilde{r}'}} &\text{if $ |t|\ge1 $.}
    \end{cases}
  \end{equation}
\end{corollary}

Analogous estimates hold for the product of $k$ hyperbolic spaces,
and more general estimates can be obtained in general
for products of Damek-Ricci spaces;
this will be the object of future work.

\begin{remark}\label{rem:strich}
  We recall that the decay rate $\sim |t|^{-\frac32}$ 
  on $\mathbb{H}^{n}$ for large times
  is sharp (see \cite{AnkerPierfelice09-a}).
  Thus we notice a new phenomenon, indeed, 
  the decay for large $t$ on
  $\mathbb{H}^{m}\times \mathbb{H}^{n}$ is faster than on
  $\mathbb{H}^{m+n}$.
  More generally, we can consider the product of $k$ real hyperbolic
  spaces ($m_{j}\ge2$)
  \begin{equation*}
    X=\mathbb{H}^{m_{1}}\times \cdots \times \mathbb{H}^{m_{k}}
  \end{equation*}
  and we obtain a decay of order 
  \begin{equation*}
    |u|\lesssim |t|^{-\frac12\sum m_{j}}
  \end{equation*}
  for small times, while the decay for large times is
  \begin{equation*}
    |u|\lesssim |t|^{-\frac32k}.
  \end{equation*}
  In particular, for spaces of the same dimension $m\ge2$
  \begin{equation*}
    X=\mathbb{H}^{m}\times \cdots \times\mathbb{H}^{m}
  \end{equation*}
  the total dimension is $mk$ but we get a decay rate 
  $\sim|t|^{-\frac32k}$ for large times. Thus if $m>3$,
  the decay rate is slower than in the euclidean case of
  the same dimension $\mathbb{R}^{mk}$, where one has
  $\sim|t|^{-\frac{mk}{2}}$. 
  On the other hand, if $m=3$ we obtain exactly the same
  decay as in the euclidean case, and if $m=2$ a better decay.
\end{remark}

As already revealed in \cite{AnkerPierfelice09-a}, the range of
exponents allowed in the dispersive estimate on $\mathbb{H}^{n}$
is wider than in the euclidean case. This reflects in a much
wider range for the Strichartz admissible indices.
Indeed, for the nonhomogeneous equation
on $X=\mathbb{H}^{m}\times \mathbb{H}^{n}$
\begin{equation}\label{eq:nonhom}
  iu_{t}-\Delta_{X}u=F(t,x,y),\qquad
  u(0,x,y)=f(x,y)
\end{equation}
we have the following result:

\begin{corollary}\label{cor:strichHnHm}
  Let $\left(\frac1p,\frac1q\right)$ and
  $\left(\frac1{\widetilde{p}},\frac1{\widetilde{q}}\right)$
  belong to the triangle
  \begin{equation}\label{eq:Tn}
    T =
    \left\{
      \left(\frac1p,\frac1q\right)\in
      \left(0,\frac12\right]\times\left(0,\frac12\right]
    \ \text{s.t. }\frac2p+\frac {m+n}q\ge \frac {m+n}2
    \right\}
    \cup 
    \left\{\left(0,\frac12\right)\right\}.
  \end{equation}
  Then the solution $u(t,x,y)$ of equation \eqref{eq:nonhom}
  on $X=\mathbb{H}^{m}\times \mathbb{H}^{n}$, $m,n\ge2$
  satisfies the estimate
  \begin{equation}\label{eq:strHnHn}
    \|u\|_{L^{p}_{t}L^{q}(X)}
    \lesssim
    \|f\|_{L^{2}(X)}
    +\|u\|_{L^{\widetilde{p}'}_{t}L^{\widetilde{q}'}(X)}.
  \end{equation}
\end{corollary}

We recall that in the euclidean case the range of admissible
indices is limited to the lower side of the triangle $T$
defined in \eqref{eq:Tn}. It is not difficult to extend 
Corollary \ref{cor:strichHnHm} to the product of $k$ real
hyperbolic spaces.

To conclude the paper we apply our estimates
to the nonlinear Schr\"odinger
equation on $X=\mathbb{H}^{m}\times \mathbb{H}^{n}$
\begin{equation}\label{eq:NLSHn}
  iu_{t}-\Delta_{X}u=F(u).
\end{equation}
We shall limit ourself here to the $L^{2}$ well posedness, but 
an analogous $H^{1}$ theory with scattering holds for suitable 
gauge invariant or defocusing type nonlinearity.
We recall that on the hyperbolic spaces $\mathbb{H}^n$, 
under the additional assumptions of radial symmetry on the data
and gauge invariance or defocusing type,
scattering properties were studied in \cite{BCS},
using the weighted radial Strichartz estimates
obtained in \cite{Ba} for $n=3$
and in \cite{P2} for $n\ge3$.  
Scattering for general power nonlinearities without
gauge invariance, and with small $L^{2}$ or $H^{1}$
data, was obtained in \cite{AnkerPierfelice09-a},  
using suitable generalized Strichartz estimates.

Here we consider a nonlinear term satisfying
\begin{equation}\label{eq:Fu}
  |F(u)|\le C|u|^{\gamma},\qquad
  |F(u)-F(v)|\le C(|u|+|v|)^{\gamma-1}|u-v|
\end{equation}
for some real $\gamma\ge1$, without gauge invariance or sign
assumptions. Then we have:

\begin{theorem}\label{the:WPNLS}
  Let $X=\mathbb{H}^{m}\times \mathbb{H}^{n}$, $m,n\ge2$.
  Assume $\gamma\le 1+ \frac4{m+n}$. Then, for all small data
  $f\in L^{2}(X)$, equation \eqref{eq:NLSHn} has a global
  unique solution, continuous with values in $L^{2}$, which
  in addition has the scattering property: there exist
  $u_{\pm}\in L^{2}$ such that
  \begin{equation}\label{eq:scat}
    \|u-e^{it \Delta_{X}}u_{\pm}\|_{L^{2}(X)}\to0 
    \quad\text{as}\qquad t\to\pm \infty.
  \end{equation}
  For large $L^{2}$ data and $\gamma<1+\frac4{m+n}$ the Cauchy
  problem is locally well posed.
\end{theorem}

\section{Proof of Theorem \ref{the:1}}\label{sec:proof_of_theorems} 

As mentioned in the Introduction, the proof of the
Theorem is completely elementary and is based on the
factorization
\begin{equation*}
  e^{itL}=e^{itH}e^{itK}
\end{equation*}
with the two flows acting on independent variables
$x\in M$ and $y\in N$. Thus we can write
\begin{equation}\label{eq:half}
  \|e^{itL}f\|_{L^{r}_{x,y}}=
  \left\|\|e^{itH}e^{itK}f\|_{L^{r}_{x}}\right\|_{L^{r}_{y}}
  \le C_{0}|t|^{-a}\|e^{itK}f\|_{L^{\widetilde{r}}_{x}L^{r}_{y}}
\end{equation}
by the first part of assumption \eqref{eq:assdis}. Now we notice
that
\begin{equation}\label{eq:LrLr}
  \|g(x,y)\|_{L^{\widetilde{r}}_{x}L^{r}_{y}}\le 
  \|g(x,y)\|_{L^{r}_{y}L^{\widetilde{r}}_{x}}\qquad
  \ \text{provided $1\le \widetilde{r}\le r\le \infty$}.
\end{equation}
Inequality \eqref{eq:LrLr} is obvious in the endpoint cases
$r=\widetilde{r}=1$ and $r=\widetilde{r}=\infty$, and in the
case $\widetilde{r}=1$, $r=\infty$ it reduces to
\begin{equation*}
  \sup_{y\in N}\int_{M}|g(x,y)|dx\le 
  \int_{M}\sup_{y\in N}|g(x,y)|dx
\end{equation*}
which is also obvious. The remaining cases follow by complex
interpolation.

Now we can continue \eqref{eq:half} using the second part of assumption
\eqref{eq:assdis} and we obtain
\begin{equation*}
  \le C_{0}|t|^{-a}\|e^{itK}f\|_{L^{r}_{y}L^{\widetilde{r}}_{x}}
  \le C_{0}|t|^{-a}\cdot C_{0}|t|^{-b}\|f\|_{L^{\widetilde{r}}_{x,y}}
\end{equation*}
and we obtain the estimate \eqref{eq:Ldis}.

\section{Sketch of the proof of Proposition \ref{pro:stri}}
\label{sec:strichartz} 

We follow the standard strategy developed by Kato, Ginibre--Velo
and Keel--Tao. 
For simplicity, we give the argument for the case $I=\mathbb{R}$
and $a+b>1$. the remaining cases are analogous.

Consider the operator
\begin{equation*}
Tf(t,x)=e^{\hspace{.25mm}i\hspace{.25mm}t\hspace{.25mm}L}f
\end{equation*}
and its formal $L^{2}$  adjoint
\begin{equation*}
T^*F=\int_{-\infty}^{+\infty}
e^{-i\hspace{.25mm}s\hspace{.25mm}L}F(s)\,ds,\qquad
F:\mathbb{R}\times X\to \mathbb{C}.
\end{equation*}
The first step of the method consists in proving
the \,$L_t^{p'}\!L_X^{q'}\!\to L_t^pL_X^q$ \,boundedness
of the operator
\begin{equation}\label{TT*}
TT^*F=\int_{-\infty}^{+\infty}
e^{\hspace{.25mm}i\hspace{.25mm}(t-s)\hspace{.25mm}L}F(s)\,ds
\end{equation}
and of its truncated version
\begin{equation}\label{truncated}
\widetilde{TT^*}F=\int_0^t
e^{\hspace{.25mm}i\hspace{.25mm}(t-s)\hspace{.25mm}L}F(s)\,ds\,,
\end{equation}
for every admissible pair $(p,q)$.
The endpoint $(\frac1p,\frac1q)\!=\!(0,\frac12)$ is settled by $L^2$ conservation
and the endpoint 
$(\frac1p,\frac1q)\!=\!(\frac12,\frac12\!-\!\frac1{2(a+b)})$
will be handled at the end.
Thus we are left with the pairs $(p,q)$
such that \,$\frac12\!-\!\frac1{2(a+b)}\!<\!\frac1q\!<\!\frac12$.
According to the dispersive estimates in Theorem \ref{the:1},
the $L_t^pL_x^q$ norms of \eqref{TT*} and \eqref{truncated}
are bounded above by
\begin{equation}\label{HLS}
\Bigl\|\,\int\hspace{-1mm}
|t\!-\!s|^{-\sigma(q)}\,
\|F(s)\|_{L_x^{q'}}\,\Bigr\|_{L_t^p}\!,\qquad
\sigma(q)=(a+b)\left(1-\frac2q\right).
\end{equation}
The convolution kernel
\,$|t\!-\!s|^{-\sigma(q)}$
\,on $\mathbb{R}$ defines a bounded operator
from $L_s^{p_1}$ to $L_t^{p_2}$,
for $p$ the first element of the admissible couple $(p,q)$,
and this proves the estimate in the non-endpoint case.
Consider eventually the endpoint
$(\frac1p,\frac1q)\!=\!(\frac12,\frac12\!-\!\frac1{2(a+b)})$;
then we can proceed exactly as in \cite{KeelTao98-a} by
splitting the time integral in dyadic regions $|t-s|\sim 2^{j}$,
$j\in \mathbb{Z}$.
Indices are finally decoupled, using the $TT^*$ argument.

\section{Proof of the corollaries}\label{sec:cor}

The proof of Corollaries \ref{cor:3}--\ref{cor:6}
is a direct application of Theorem \ref{the:1}, combined with
dispersive estimates from different papers. More precisely:

\begin{enumerate}
  \item For Corollary \ref{cor:3}, we use $n$ times
  the dispersive estimate
  \begin{equation*}
    |P_{ac}e^{itH}f|\lesssim|t|^{-1/2}\|f\|_{L^{1}(\mathbb{R})}
  \end{equation*}
  where $H$ is the Schr\"odinger operator on $\mathbb{R}$
  \begin{equation*}
    H=-\frac{d^{2}}{dx^{2}}+V(x),\qquad
    (1+|x|)^{2}V\in L^{1}(\mathbb{R})
  \end{equation*}
  and $P_{ac}$ is the projection on the absolutely continuous
  space associated to $H$ (see \cite{GoldbergSchlag04-b}
  and \cite{DanconaFanelli06-a}). From the general theory
  it is known that all eigenvalues (if present) must be
  nonnegative. The additional assumption
  $V\ge0$ ensures that no eigenvalues exist so that
  the projection $P_{ac}$ is not necessary.
  
  \item For Corollaries \ref{cor:4} and
  \ref{cor:5}, we use the results of 
  \cite{Yajima95-a}, where it is proved that under assumption
  \eqref{eq:yaj} the wave operator associated to $H=-\Delta+V$
  is bounded on $L^{p}$. In particular, this gives dispersive
  estimates for the Schr\"odinger equation of the form
  \begin{equation*}
    |e^{itH}f|\lesssim|t|^{-n/2}\|f\|_{L^{1}}.
  \end{equation*}
  To our knowledge,
  Yajima's conditions are the best known for large space dimension
  $n\ge4$.

  \item Corollary \ref{cor:6} follows easily from \eqref{eq:AP}.
  Notice that the wider range of exponents compared
  with the euclidean case
  is due to the Kunze-Stein phenomenon
  \begin{equation*}
    \|f*g\|_{L^{q',\infty}(\mathbb{H}^{n})}\lesssim
    \|f\|_{L^{q'}(\mathbb{H}^{n})}
    \|g\|_{L^{q'}(\mathbb{H}^{n})}\qquad
    \forall q>2
  \end{equation*}
  (see \cite{Cowling97-a}, \cite{Ionescu00-a}); however we need here
  only the endpoint case $q=\widetilde{q}=\infty$ in order
  to verify the assumptions of Theorem \ref{the:1}.

  \item Corollary \ref{cor:strichHnHm} is proved by the same 
  $TT^{*}$ method as sketched in Section \ref{sec:strichartz},
  however with an important difference since now the rate of
  decay for small and large times is different. As above,
  we consider the operators
  \begin{equation}\label{TT*}
    TT^*F(t,x,y)=\int_{-\infty}^{+\infty}
    e^{\hspace{.25mm}i\hspace{.25mm}(t-s)
    \hspace{.25mm}{\Delta}_X}F(s,x,y)\,ds
  \end{equation}
  and
  \begin{equation}\label{truncated}
    \widetilde{TT^*}F(t,x,y)=\int_0^t
    e^{\hspace{.25mm}i\hspace{.25mm}(t-s)
    \hspace{.25mm}{\Delta}_X}F(s,x,y)\,ds\,,
  \end{equation}
  for every admissible pair $(p,q)$.
  The endpoint $(\frac1p,\frac1q)\!=\!(0,\frac12)$ is true by 
  $L^2$ conservation.
  The endpoint $(\frac1p,\frac1q)\!=\!(\frac12,\frac12\!-\!\frac1{m+n})$
  for $m+n\!\ge\!3$ is handled by the standard method of
  \cite{KeelTao98-a} applied to the truncated $\widetilde{TT^{*}}$
  directly.
  Finally consider the pairs $(p,q)$
  such that $\frac12\!-\!\frac1{m+n}\!<\!\frac1q\!<\!\frac12$
  and $(\frac12\!-\!\frac1q)\frac {m+n}2\!\le\!\frac1p\!\le\!\frac12$,
  for which it is sufficient to study $TT^{*}$.
  According to the dispersive estimates in Corollary \ref{cor:6},
  the $L_t^pL^q(X)$ norms of \eqref{TT*} and \eqref{truncated}
  are bounded above by
  \begin{equation}\label{HLS}
  \Bigl\|\,\int_{\,|t-s|\ge1}\hspace{-1mm}
  |t\!-\!s|^{-3}\,
  \|F(s)\|_{L_x^{q'}}\,\Bigr\|_{L_t^p}\!
  +\,\Bigl\|\,\int_{\,|t-s|\le1}\hspace{-1mm}
  |t\!-\!s|^{-(\frac12-\frac1q)(m+n)}\,
  \|F(s)\|_{L_x^{q'}}\,\Bigr\|_{L_t^p}\,.
  \end{equation}
  The convolution kernel in the first integral
  $|t\!-\!s|^{-3}\,1\hspace{-1mm}\text{l}_{\,\{|t-s|\ge1\}}$
  on $\mathbb{R}$ defines a bounded operator
  from $L_s^{p_1}$ to $L_t^{p_2}$,
  for all $1\!\le\!p_1\!\le\!p_2\!\le\!\infty$,
  in particular from $L_s^{p'}$ to $L_t^p$,
  for all $2\!\le\!p\!\le\!\infty$.
  The convolution kernel in the second integral
  $|t\!-\!s|^{-(\frac12-\frac1q)\hspace{.25mm}(m+n)}\,
  1\hspace{-1mm}\text{l}_{\,\{|t-s|\le1\}}$
  defines a bounded operator
  from $L_s^{p_1}$ to $L_t^{p_2}$,
  for all $1\!<\!p_1,p_2\!<\!\infty$
  such that $0\!\le\!\frac1{p_1}\!-\!\frac1{p_2}\!
  \le\!1\!-\!(\frac12\!-\!\frac1q)(m+n)$,
  in particular from $L_s^{p'}$ to $L_t^p$,
  for all $2\!\le\!p\!<\!\infty$ such that
  $\frac1p\!\ge\!(\frac12\!-\!\frac1q)\frac {m+n}2$.
  This proves the result for all dual estimates with
  $(p,q)=(\widetilde{p},\widetilde{q})$. The standard $TT^{*}$
  argument allows to decouple the pairs and conclude the proof.
\end{enumerate}

\section{The nonlinear Schr\"odinger equation on $\mathbb{H}^{m}\times \mathbb{H}^{n}$}\label{sec:NLS} 

We follows the same strategy used in \cite{AnkerPierfelice09-a}.
We resume the standard fixed point method based on Strichartz estimates.
Define \,$u\!=\!\Phi(v)$ \,as the solution to the Cauchy problem
\begin{equation}\label{Phi}
\begin{cases}
\;i\,\partial_tu(t,x,y)+{\Delta}_X u(t,x,y)=F(v(t,x,y))\,,\\
\;u(0,x,y)=f(x,y)\,,\\
\end{cases}
\end{equation}
which is given by Duhamel's formula:
\begin{equation*}
u(t,x,y)=\,
e^{\hspace{.25mm}i\hspace{.25mm}t\hspace{.25mm}{\Delta}_X}\hspace{-.25mm}f(x,y)
\hspace{.25mm}+\int_0^t\!
e^{\hspace{.25mm}i\hspace{.25mm}(t-s)\hspace{.25mm}{\Delta}_X}
F(v(s,x,y))\,ds\,.
\end{equation*}
According to  \eqref{eq:strHnHn},
we have the following Strichartz estimate
\vspace{1mm}
\begin{equation}\label{StrichartzL2v1}
\|u\|_{L_t^\infty L^2(X)\vphantom{L_t^{\tilde p'}}}\!
+\|u\|_{L_t^pL^q(X)\vphantom{L_t^{\tilde p'}}}
\le C\,\|f\|_{L^2(X)\vphantom{L_t^{\tilde p'}}}\!
+C\,\|F(v)\|_{L_t^{\tilde p'}\!L^{\tilde q '}(X)}
\end{equation}
for all $(\frac1p,\frac1q)$
and $(\frac1{\tilde p},\frac1{\tilde q})$
in the triangle $T$,
which amounts to the conditions
\begin{equation}\label{admissibilityL2}
\begin{cases}
\;2\!\le\!p,q\!\le\!\infty\text{ \;such that \;}
\frac\beta p\!+\frac{m+n}{q}=\!\frac{m+n}{2}
\text{ \;for some \;}0\!<\!\beta\!\le\!2\,,\\
\;2\!\le\!\tilde p,\tilde q\!\le\!\infty\text{ \;such that \;}
\frac{\widetilde{\beta}}{\widetilde{p}}\!+\frac{m+n}{\widetilde{q}}=\!\frac{m+n}{2}
\text{ \;for some \;}0\!<\!\tilde\beta\!\le\!2\,.\\
\end{cases}
\end{equation}
Moreover
\begin{equation*}
\|F(v)\|_{L_t^{\tilde p'}\!L^{\tilde q '}(X)}
\le C\,\|\,|v|^\gamma\|_{L_t^{\tilde p'}\!L^{\tilde q '}(X)}
\le C\,\|v\|_{
L_{\,t}^{\tilde p'\gamma}\!L^{\tilde q '\gamma}(X)
}^{\,\gamma}
\end{equation*}
by our nonlinear assumption.
Thus
\vspace{1mm}
\begin{equation}\label{StrichartzL2v2}
\|u\|_{L_t^\infty L^2(X)\vphantom{L_t^{\tilde p'}}}\!
+\|u\|_{L_t^pL^q(X)\vphantom{L_t^{\tilde p'}}}
\le C\,\|f\|_{L^2(X)\vphantom{L_t^{\tilde p'}}}\!
+C\,\|v\|_{
L_{\,t}^{\tilde p'\gamma}\!L^{\tilde q '\gamma}(X)
}^{\,\gamma}\,.
\end{equation}
In order to remain within the same function space,
we require in addition
\begin{equation}
\label{selfmappingL2}
p=\tilde p' \gamma,	
\;q=\tilde q' \gamma.
\end{equation}
It is easily checked that all these conditions are fulfilled
if we take for instance
\begin{equation*}\textstyle
0<\beta=\tilde\beta\le2
\quad\text{such that}\quad
\gamma=1\!+\!\frac{2\,\beta}{m+n}
\quad\text{and}\quad
p=q=\tilde p=\tilde q=1\!+\!\gamma=2\!+\!\frac{2\,\beta}{m+n}\,.
\end{equation*}
For such a choice, $\Phi$ maps
$L^\infty(\mathbb{R};L^2(X))
\cap L^p(\mathbb{R};L^q(X))$
into itself, and actually
\linebreak
$Y\!=C(\mathbb{R};L^2(X))
\cap L^p(\mathbb{R};L^q(X))$
into itself.
Since $Y$ is a Banach space for the norm
\begin{equation*}
\|u\|_Y=\|u\|_{L_t^\infty L^2}+\|u\|_{L_t^pL^q}\,,
\end{equation*}
it remains for us to show that $\Phi$ is a contraction in the ball
\begin{equation*}
Y_\varepsilon=\{\,u\!\in\!Y\mid\|u\|_Y\!\le\!\varepsilon\,\}\,,
\end{equation*}
provided $\varepsilon\!>\!0$ and \,$\|f\|_{L^2}$ are sufficiently small.
Let $v,\tilde v\!\in\!X$ and $u\!=\!\Phi(v)$, $\tilde u\!=\!\Phi(\tilde v)$.
Arguying as above and using in addition H\"older's inequality,
we estimate
\vspace{1mm}
\begin{equation*}
\begin{aligned}
\|u-\tilde u\|_{Y\vphantom{L_t^{\tilde p'}}}\,
&\le\,C\;\|F(v)-F(\tilde v)\|_{L_t^{\tilde p'}\!L^{\tilde q'}}\\
&\le\,C\;\|
\{\,|v|^{\gamma-1}\!+|\tilde v|^{\gamma-1}\}\,|v\!-\!\tilde v|
\|_{L_t^{\tilde p'}\!L^{\tilde q'}}\\
&\le\,C\;
\bigl\{\|v\|_{L_t^pL^q}^{\,\gamma-1}\!
+\|\tilde v\|_{L_t^pL^q}^{\,\gamma-1}\bigr\}\,\|v-\tilde v\|_{L_t^pL^q\vphantom{L_t^{\tilde p'}}},
\end{aligned}
\end{equation*}
hence
\begin{equation}\label{contractionL2}
\|u-\tilde u\|_Y
\le C\,\bigl(\,\|v\|_Y^{\gamma-1}\!+\|\tilde v\|_Y^{\gamma-1}\,\bigr)\,
\|v-\tilde v\|_{Y}\,.
\end{equation}
If we assume $\|v\|_Y\!\le\!\varepsilon$,
$\|\tilde v\|_Y\!\le\!\varepsilon$
and $\|f\|_{L^2}\!\le\!\delta$,
then \eqref{StrichartzL2v2} and \eqref{contractionL2} yield
\begin{equation*}
\|u\|_Y\le C\,\delta+C\,\varepsilon^\gamma\,,\;
\|\tilde u\|_Y\le C\,\delta+C\,\varepsilon^\gamma
\quad\text{and}\quad
\|u-\tilde u\|_Y\le 2\;C\,\varepsilon^{\gamma-1}\,\|v-\tilde v\|_Y\,.
\end{equation*}
Thus
\begin{equation*}\textstyle
\|u\|_Y\le\varepsilon\,,\;
\|\tilde u\|_Y\le\varepsilon
\quad\text{and}\quad
\|\,u-\tilde u\,\|_Y\le\frac12\,\|\,v-\tilde v\,\|_Y
\end{equation*}
if \,$C\,\varepsilon^{\gamma-1}\!\le\frac14$
\,and \,$C\,\delta\le\frac34\,\varepsilon$\,.
We conclude by applying the fixed point theorem
in the complete metric space $Y_\varepsilon$.
Hence,
for \,$1\!<\!\gamma\!\le\!1\!+\!\frac4{m+n}$
\,and small $L^2$ data,
the Cauchy problem (\ref{eq:Fu})
has a unique solution $u(t,x,y)$
in $C(\mathbb{R};L^2(X))\cap
L^p(\mathbb{R};L^q(X))$,
for the above choice of a suitable pair $(p,q)$.
Scattering will follow from the Cauchy criterion\,:
\begin{quote}
If \,$\|z(t_1)-z(t_2)\|_{L^2(X)}\to 0$
\,as \,$t_1,t_2\!\to\!+\infty$\,,
then there exists $z_+\!\in\!L^2$
such that \,$\|z(t)-z_+\|_{L^2(X)}\to 0$
\,as \,$t\!\to\!+\infty$\,.
\end{quote}
In our case $z(t,x,y)=e^{-i\hspace{.25mm}t\hspace{.25mm}\Delta_{X}}u(t,x,y)$.
So if we prove that
\begin{equation*}
\|e^{-i\hspace{.25mm}t_2\hspace{.25mm}\Delta_{X}}u(t_2)
-e^{-i\hspace{.25mm}t_1\hspace{.25mm}\Delta_{X}}u(t_1)\|_{L^2(X)}
\to0
\quad\text{as}\quad
t_1\le t_2\to\pm\infty\,,
\end{equation*}
we can conclude that the global solution $u(t,x))$
has the scattering property stated above.
Using our Strichartz estimates \eqref{eq:strHnHn}, we get
\begin{equation*}\begin{aligned}
\bigl\|e^{-i\hspace{.25mm}t_2\hspace{.25mm}\Delta_{X}}u(t_2)
-e^{-i\hspace{.25mm}t_1\hspace{.25mm}\Delta_{X}}u(t_1)\bigr\|_{L^2{(X)}}
&=\,\Bigl\|\,\int_{t_1}^{t_2}\hspace{-1mm}
e^{-i\hspace{.25mm}s\hspace{.25mm}\Delta_{X}}
F(u(s))\,ds\,\Bigr\|_{L^2{(X)}}\\
&\le\,\bigl\|u\bigr\|_{L^p([t_1,t_2];L^q(X))}^{\,\gamma}\,.
\end{aligned}\end{equation*}
Since $u(t,x,y)\!\in\!L^p(\mathbb{R};L^q(X))$,
the last expression vanishes
as \,$t_1\!\le\!t_2$ \,tend both to $+\infty$ or $-\infty$\,.

In the subcritical case $\gamma\!<\!1\!+\!\frac{4}{m+n}$,
one can prove in a similar way
local well--posedness in $L^2$ for arbitrary data $f$.
Specifically,
we restrict to a small time interval $I\!=\![-T,+T]$
and proceed as above,
except that we increase $\tilde\beta\!\in\!(\beta,2\,]$
and ${\tilde p}\!=\!\frac{\tilde\beta}{\beta}\,p$ accordingly,
and that we apply in addition H\"older's inequality in time.


\begin{thebibliography}{10}

\bibitem{AnkerPierfelice09-a}
Jean-Philippe Anker and Vittoria Pierfelice.
\newblock Nonlinear {S}chr{\"o}dinger equation on real hyperbolic spaces.
\newblock {\em Ann. Inst. H. Poincar\'e Anal. Non Lin\'eaire},
  26(5):1853--1869, 2009.
  
  
  \bibitem{Ba}
Valeria Banica\,:
\textit{The nonlinear Schr\"odinger equation
on the hyperbolic space\/},
Comm. P.D.E. 32 (2007), no. 10, 1643--1677
[arXiv:math/0406058]


\bibitem{BCS}
Valeria  Banica, R. Carles, G. Staffilani\,:
\textit{Scattering theory for radial nonlinear Schr\"odinger equations
on hyperbolic space\/},
Geom. Funct. Anal. (to appear) [arXiv:math/0607186]

\bibitem{Cowling97-a}
Michael Cowling.
\newblock Herz's ``principe de majoration'' and the {K}unze-{S}tein phenomenon.
\newblock In {\em Harmonic analysis and number theory ({M}ontreal, {PQ},
  1996)}, volume~21 of {\em CMS Conf. Proc.}, pages 73--88. Amer. Math. Soc.,
  Providence, RI, 1997.

\bibitem{DanconaFanelli06-a}
Piero D'Ancona and Luca Fanelli.
\newblock {$L^p$}-boundedness of the wave operator for the one dimensional
  {S}chr{\"o}dinger operator.
\newblock {\em Comm. Math. Phys.}, 268(2):415--438, 2006.

\bibitem{GinibreVelo95-a}
Jean Ginibre and Giorgio Velo.
\newblock Generalized {S}trichartz inequalities for the wave equation.
\newblock In {\em Partial differential operators and mathematical physics
  (Holzhau, 1994)}, volume~78 of {\em Oper. Theory Adv. Appl.}, pages 153--160.
  Birkh\"{a}user, Basel, 1995.

\bibitem{GoldbergSchlag04-b}
Michael Goldberg and Wilhelm Schlag.
\newblock Dispersive estimates for {S}chr\"{o}dinger operators in dimensions
  one and three.
\newblock {\em Comm. Math. Phys.}, 251(1):157--178, 2004.

\bibitem{Ionescu00-a}
Alexandru~D. Ionescu.
\newblock An endpoint estimate for the {K}unze-{S}tein phenomenon and related
  maximal operators.
\newblock {\em Ann. of Math. (2)}, 152(1):259--275, 2000.

\bibitem{Kato87-a}
Tosio Kato.
\newblock On nonlinear {S}chr\"{o}dinger equations.
\newblock {\em Ann. Inst. H. Poincar\'{e} Sect. A (N.S.)}, 46(1):113--129,
  1987.

\bibitem{KeelTao98-a}
Markus Keel and Terence Tao.
\newblock Endpoint {S}trichartz estimates.
\newblock {\em Amer. J. Math.}, 120(5):955--980, 1998.

\bibitem{Pierfelice06-b}
Vittoria Pierfelice.
\newblock Weighted {S}trichartz estimates for the radial perturbed
  {S}chr{\"o}dinger equation on the hyperbolic space.
\newblock {\em Manuscripta Math.}, 120(4):377--389, 2006.



\bibitem{P2}
Vittoria Pierfelice\,:
\textit{Weighted Strichartz estimates
for the Schr\"odinger and wave equations
on Damek--Ricci spaces\/},
Math. Z. 260 (2008) , no. 2, 377-392

\bibitem{Yajima95-a}
Kenji Yajima.
\newblock The {$W\sp {k,p}$}-continuity of wave operators for schr\"{o}dinger
  operators. iii. even-dimensional cases $m\geq 4$.
\newblock {\em J. Math. Sci. Univ. Tokyo}, 2(2):311--346, 1995.

\bibitem{Yajima99-a}
Kenji Yajima.
\newblock {$L\sp p$}-boundedness of wave operators for two-dimensional
  {S}chr\"{o}dinger operators.
\newblock {\em Comm. Math. Phys.}, 208(1):125--152, 1999.

\end{thebibliography}
\end{document}